\documentclass[11pt,
]{article}

\usepackage{graphicx,color}
\usepackage{times}

\usepackage[colorlinks=true, breaklinks=true, pdfstartview=FitV, linkcolor=red, citecolor=blue, urlcolor=black]
{hyperref}

\usepackage[T1]{fontenc}

\begin{document}

\title{ Brouwer's satellite solution redux\thanks{%
Submitted to Celestial Mechanics and Dynamical Astronomy}
}

\author{Martin Lara%
\thanks{GRUCACI, University of La Rioja, and Space Dynamics Group -- UPM}
\thanks{\tt{mlara0@gmail.com} }
}

\date{\today }

\maketitle

\begin{abstract}
Brouwer's solution to the artificial satellite problem is revisited to show that the complete Hamiltonian reduction is rather achieved in the plain Poin\-ca\-r\'e's style, through a single canonical transformation, than using a sequence of partial reductions based on von Zeipel's alternative for dealing with perturbed degenerate Hamiltonian systems.
\end{abstract}



\section{Introduction}

Brouwer's analytical solution to the artificial satellite problem \cite{Brouwer1959} based on von Zeipel's partial reduction method for dealing with perturbed degenerate Hamiltonians \cite{vonZeipel1916} fiercely resists obsolescence sixty years after publication. Indeed, in spite of the spectacular increase of the computational power, widespread software packages for approximate ephemeris prediction still rely on Brouwer's seminal results \cite{HootsRoehrich1980,ValladoCrawfordHujsakKelso2006}. Furthermore, the success of Brouwer's closed-form solution among practitioners as well as the reputation gained among theorists by Brouwer's stepwise normalization, make that these days some authors designate the method as ``the von Zeipel-Brouwer theory'' \cite{FerrazMello2007}. 
\par

Since then, the merits of Brouwer's decomposition of the solution of perturbed Keplerian motion into secular, long-, and short-period effects, seem not to have been questioned. Moreover, after the invention of Hamiltonian simplification methods \cite{Deprit1981}, it has been suggested that carrying out additional decompositions, thus increasing the number of canonical transformations, could be the proper way to success in the search for separable perturbation Hamiltonians of celestial mechanics problems \cite{DepritMiller1989}. Conversely, it has been recently pointed out that the use of Hamiltonian simplification procedures could be merely optional in the construction of higher-order analytical approximations to the satellite problem \cite{Lara2020}. Then, it emerges the question of which is the real value of splitting a normalization procedure, either partial or complete, into several different stages, a topic that may well deserve additional study.
\par

We walk a step in that direction and, in the classical style of Poincar\'e's ``new method'' \cite{Poincare1892vII}, undertake the construction of Brouwer's second-order completely reduced Hamiltonian of the main problem in the artificial satellite theory (the $J_2$ problem) by means of a single canonical transformation. The difficulties stemming from the degeneracy of the Kepler Hamiltonian, who has two null frequencies, are easily overcome with the addition of suitable integration constants to the generating function of the transformation that yields the complete Hamiltonian reduction.
\par

The use of arbitrary functions in the construction of perturbation solutions is not new at all \cite{Morrison1966}. In fact, it can be traced back to Poincar\'e's efforts in approaching degenerate perturbation problems \cite[Chap.~XI]{Poincare1892vII}. They also play a fundamental role in a \emph{reverse} partial normalization process in which the angular momentum is normalized in first place \cite{AlfriendCoffey1984,LaraSanJuanLopezOchoa2013c,SanJuanOrtigosaLopezOchoaLopez2013,Lara2020}. On the other hand, in spite of average perturbation Hamiltonians do not exist in general \cite{FerrazMello1999}, the use of arbitrary constants to guarantee the purely periodic nature of the generating function became customary in attempts to bring the mean elements dynamics as close as possible to the true average dynamics \cite{MetrisExertier1995,Steichen1998,LaraSanJuanFolcikCefola2011}.
\par

To the first order, the construction of Brouwer's closed-form solution  by means a single transformation amounts to the sum of the two transformations computed by Brouwer for the short- and long-period elimination. This is due to the linearization provided by the first order truncation of a perturbation theory. In view of no differences arise between Brouwer's and the current approach when the periodic corrections are constrained to first order effects, we feel compelled to supplement Brouwer's analytical solution with second-order periodic corrections, yet limited to the $J_2$ contributions. We compare our results with corresponding corrections obtained in the traditional way, in which the normalization is split into the elimination of the parallax, the elimination of the perigee, and the Delaunay normalization \cite{AlfriendCoffey1984,Lara2019UR}. At the second order the single transformation is no longer the addition of the different canonical transformations. As expected, the (single) periodic corrections are now much more involved than those corresponding to each of the partial reductions or simplifications, and are also more intricate than the composition of all of them. However, the length of the series defining the solution is only a part of the whole picture, and we found clear computational advantages in the evaluation of the single transformation. The improvements stem from the fact that many inclination polynomials pertaining to the periodic corrections admit factorization. Because common factors repeat many times throughout the corrections, the compiler is able to perform a higher optimization of the code in the case of the single transformation than in the case in which the transformation is split in different stages. 
\par

In the process of computing the second-order corrections, we will recognize how artificial the controversy created about the integration of the equation of the center was. Indeed, difficulties confronted by researchers involved in the automatization of celestial mechanics computations were, in fact, derived from their own programming strategies \cite{DepritRom1970,Jefferys1971}. On the contrary, the trouble had been easily sidestepped by celestial mechanicians relying on traditional hand computations \cite{Kozai1962,Aksnes1971}. We will show that standard integration by parts reduces the equation of the center issue to the well known integration of cosine functions in elliptic motion \cite{Kozai1962AJ,Tisserand1889}. We hasten to say that the controversy was in no way futile since it provoked the appearance of Hamiltonian simplification methods and led to the development of sophisticated computational strategies \cite{Healy2000}.\footnote{A brief review of the history of Hamiltonian simplification methods can be found in the Introduction of \cite{Lara2019CMDA}.}
\par

In order to fully determine the second-order term of the generating function of Brouwer's theory, the third-order term of the completely reduced Hamiltonian needs to be previously specified. The use of higher-order secular terms should improve further the long-term performance of Brouwer's solution. However, to be effective in the propagation of an initial state vector, the initialization constants of the analytical solution, and, in particular, the secular frequencies, must be computed within comparable accuracy to that of the secular terms. Rather than carrying out the long and tedious computations required in the determination of the third-order generating function, we take the clever shortcut proposed by Breakwell and Vagners \cite{BreakwellVagners1970}. That is, we limit the computation of third-order corrections to the case of the secular mean motion, which, besides, is directly obtained from the secular Hamiltonian. With this effortless procedure the addition of third-order secular terms clearly improves the performance of Brouwer's solution. 
\par

\section{Brouwer's complete reduction at once} \label{}

Constraining the dynamical model of the artificial satellite problem to the $J_2$ perturbation (main problem), Brouwer's gravitational Hamiltonian takes the form \cite{Brouwer1959}
\begin{equation} \label{Hoscu}
\mathcal{H}=-\frac{\mu}{2a}+\frac{\mu}{r}\frac{R_\oplus^2}{r^2}\frac{1}{2}C_{2,0}\left[1-\frac{3}{2}s^2+\frac{3}{2}s^2\cos(2f+2\omega)\right],
\end{equation}
where the Earth's gravitational field is materialized by the physical constants $\mu$, the gravitational parameter, $R_\oplus$, the equatorial radius, and $C_{2,0}=-J_2$, the non-dimensional oblateness coefficient.\footnote{Note that $k_2=-\frac{1}{2}C_{2,0}R_\oplus^2$ in Brouwer's notation.} The symbols $a$, $r$, $f$ and $\omega$, stand for semi-major axis, radius from the Earth's center of mass,  true anomaly, and argument of the perigee, respectively, whereas $s\equiv\sin{I}$ abbreviates the sine of the inclination $I$. Since we are dealing with Hamiltonian mechanics, these symbols must be understood as \emph{functions} of some set of canonical variables. In particular we assume, with Brouwer, that the Hamiltonian is written in terms of the Delaunay coordinates $\ell$, the mean anomaly, $g=\omega$, and $h$, the longitude of the ascending node, and their conjugate momenta $L=\sqrt{\mu{a}}$, $G=L(1-e^2)^{1/2}$, with $e$ denoting the eccentricity, and $H=G\cos{I}$, standing for the Delaunay action, the total angular momentum, and the projection of the angular momentum vector along the Earth's rotation axis, respectively. That $H$ is an integral of Eq.~(\ref{Hoscu}) becomes evident from the cyclic character of $h$. Besides, the Hamiltonian itself is constant because the time does not appear explicitly on it.
\par

The small value of the Earth's $J_2$ coefficient identifies Eq.~(\ref{Hoscu}) like a case of perturbed Keplerian motion, which, therefore, can be reduced to a separable Hamiltonian by perturbation methods. This is achieved by finding a canonical transformation $\mathcal{T}:(\ell,g,h,L,G,H,\epsilon)\mapsto(\ell',g',h',L',G',H')$, from osculating to mean variables, depending on the small parameter $\epsilon\sim{J}_2$, such that the transformed Hamiltonian in mean (prime) variables becomes a function of only the momenta, namely $\mathcal{H}\circ\mathcal{T}=\mathcal{K}(-,-,-,L',G',H';\epsilon)$. The transformation $\mathcal{T}$, we learned from Poincar\'e \cite{Poincare1892vII}, is derived from a determining function that is solved in the form of a Taylor series up to some truncation order of the small parameter $\epsilon$.
\par

Brouwer, for his part, after introducing the method of solution, seems to refuse approching the direct computation of the transformation $\mathcal{T}$ since the beginning, by simply declaring that
\begin{quote}
``[\dots] it is more convenient to choose a determining function in such a manner that the mean anomaly is not present in the transformed Hamiltonian while the argument of the perigee is permitted to appear.'' 
\end{quote}
Next, after invoking von Zeipel, he proceeded stepwise by partial reduction, first computing a canonical transformation that only removes the short-period terms from the Hamiltonian, and then carrying out a second canonical transformation that removes the long-period terms. In this way Brouwer outstandingly achieves the complete Hamiltonian reduction in closed form.
\par

Conversely, we ignore the presumed convenience of Brouwer's procedure and approach the perturbation problem searching for a single determining function in the original style of Poincar\'e, yet we better rely on the equivalent but more functional method of Lie transforms \cite{Hori1966,Deprit1969,DepritDeprit1999}. Thus, we write Eq.~(\ref{Hoscu}) in the usual form of a perturbation Hamiltonian $\mathcal{H}=\sum_{m\ge0}(\epsilon^m/m!)\mathcal{H}_{m,0}$, with
\[
\mathcal{H}_{0,0}=-\frac{\mu}{2a}, \qquad
\mathcal{H}_{1,0}=-\frac{\mu}{r}\frac{R_\oplus^2}{r^2}\frac{1}{2}\left[1-\frac{3}{2}s^2+\frac{3}{2}s^2\cos(2f+2\omega)\right],
\]
$\mathcal{H}_{m,0}=0$ for $m\ge2$, and $\epsilon\equiv{J}_2=-C_{2,0}$. Recall, that all the symbols are functions of the Delaunay canonical variables, in which the Lie operator $\mathcal{L}_0=\{\quad;\mathcal{H}_{0,0}\}$, where the curly brackets stand for the Poisson brackets operator, takes the simple form $\mathcal{L}_0=n\partial/\partial\ell$, where $n=\mu^2/L^3$ is the mean motion. This allows us to compute the determining function $\mathcal{W}=\sum_{m\ge0}(\epsilon^m/m!)\mathcal{W}_{m+1}$ from the sequence given by the homological equation
\begin{equation} \label{homological}
\mathcal{W}_m=\frac{1}{n}\int(\widetilde{\mathcal{H}}_{0,m}-\mathcal{H}_{0,m})\mathrm{d}\ell+\mathcal{C}_m.
\end{equation}
At each step $m$, terms $\widetilde{\mathcal{H}}_{0,m}$ in Eq.~(\ref{homological}) are known, coming either from the original Hamiltonian or stemming from intermediate computations at previous orders. Terms $\mathcal{H}_{0,m}$ are selected in such a way that they cancel those terms of $\widetilde{\mathcal{H}}_{0,m}$ pertaining to the kernel of the Lie operator. Finally, the integration ``constants'' $\mathcal{C}_m$ ---arbitrary functions of the Delaunay variables fulfilling the condition $\partial\mathcal{C}_m/\partial\ell=0$--- will be chosen like such trigonometric functions of $g$ that they prevent the appearance of purely long-period terms at the next order of the perturbation approach, in this way making feasible the complete normalization at once. The method is standard these days, and the required details can be found in textbooks as, for instance, \cite{MeyerHallOffin2009,BoccalettiPucacco1998v2}.
\par

In preparation of the solution, the equivalence
\begin{equation} \label{runity}
\frac{1}{r^j}=\frac{1}{r^2}\left(\frac{1+e\cos{f}}{p}\right)^{j-2},
\end{equation}
where $p=a\eta^2$ is the orbit parameter and $\eta=(1-e^2)^{1/2}$ is the so-called eccentricity function, is applied to the instances $j>2$ in $\widetilde{\mathcal{H}}_{0,1}=\mathcal{H}_{1,0}$, which is then written in the convenient form 
\begin{equation} \label{Ktilde1}
\widetilde{\mathcal{H}}_{0,1}=\mathcal{H}_{0,0}\frac{R_\oplus^2}{r^2}\frac{1}{\eta^2}
\sum _{i=0}^1B_i(s)\sum _{j=i}^{2i+1}(2-j^\star)^ie^{|j-2i|}\cos(jf+2ig),
\end{equation}
where $B_0=1-\frac{3}{2}s^2$, $B_1=\frac{3}{4}s^2$, and we abbreviate $j^\star\equiv{j}\bmod2$. On account of 
$j\ge{k}$ in Eq.~(\ref{Ktilde1}), we immediately verify that $\widetilde{\mathcal{H}}_{0,1}$ is not affected of purely long-period terms. Then, the complete reduction is achieved at the first order by choosing the new Hamiltonian term $\mathcal{H}_{0,1}$ like the average of $\mathcal{H}_{1,0}$ over the mean anomaly.
\par

The average is obtained in closed form with the help of the Keplerian differential relation between the true and mean anomalies $\eta{a}^2\mathrm{d}\ell=r^2\mathrm{d}f$. It is equivalent to removing all the terms with $j>0$ from Eq.~(\ref{Ktilde1}) after multiplied by the factor $r^2/(a^2\eta)$. We trivially obtain the usual result
\[
\mathcal{H}_{0,1}=\langle\widetilde{\mathcal{H}}_{0,1}\rangle_\ell
\equiv\mathcal{H}_{0,0}\frac{R_\oplus^2}{p^2}\eta\left(1-\frac{3}{2}s^2\right),
\]
which is, of course, the same expression obtained by Brouwer. Then, Eq.~(\ref{homological}) is rearranged in the form
\begin{equation} \label{homologicalf}
\mathcal{W}_1=\frac{1}{n}\left[\mathcal{H}_{0,1}\phi+\int\left(\widetilde{\mathcal{H}}_{0,1}\frac{r^2}{a^2\eta}-\mathcal{H}_{0,1}\right)\mathrm{d}f\right]+\mathcal{C}_1,
\end{equation}
where $\phi=f-\ell$ denotes the equation of the center, and the integrand in Eq.~(\ref{homologicalf}) only embraces periodic functions of $f$.  We obtain
\begin{equation} \label{W1}
\mathcal{W}_1=-G\frac{R_\oplus^2}{p^2}\frac{1}{2}\bigg[B_0\phi+
\sum _{i=0}^1B_i\sum_{j=\max(i,1)}^{2i+1}\frac{(2-j^\star)^i}{j}e^{|j-2i|}\sin(jf+2ig)\bigg]+\mathcal{C}_1,
\end{equation}
where the first term of the right hand member is the same as Brouwer's first order determining function of the short-period elimination, and $\mathcal{C}_1$ is an integration constant that is left undetermined by the time being.
\par

On account of $\mathcal{H}_{2,0}\equiv0$, the known terms at the second order of the Lie transforms approach are $\widetilde{\mathcal{H}}_{0,2}=\{\mathcal{H}_{1,0};\mathcal{W}_1\}+\{\mathcal{H}_{0,1};\mathcal{W}_1\}$, from which the terms pertaining to the kernel of the Lie operator must be cancelled by the adequate selection of $\mathcal{H}_{0,2}$. The usual choice is $\mathcal{H}_{0,2}=\langle\widetilde{\mathcal{H}}_{0,2}\rangle_\ell$, yet additional terms could be left in the new Hamiltonian in particular cases \cite{Deprit1981,Lara2019CMDA}. However, this process would leave purely long-period terms in the new Hamiltonian in addition to the secular terms, both certainly pertaining to the kernel of the Lie operator. Since this is against the total normalization criterion, purely long-period terms should vanish identically in $\widetilde{\mathcal{H}}_{0,2}$, a requirement that is achieved with the proper selection of $\mathcal{C}_1$, whose partial derivatives with respect to $g$, $G$, and $L$, appear formally in $\widetilde{\mathcal{H}}_{0,2}$.
\par

We attack the computation of the second-order of the perturbation theory by parts. To that effect, we make $\widetilde{\mathcal{H}}_{0,2}=\widetilde{\mathcal{H}}'_{0,2}+\widetilde{\mathcal{H}}^*_{0,2}$ with
$\widetilde{\mathcal{H}}'_{0,2}=\{\mathcal{H}_{1,0}+\mathcal{H}_{0,1};\mathcal{V}_1\}$ and $\widetilde{\mathcal{H}}^*_{0,2}=\{\mathcal{H}_{1,0}+\mathcal{H}_{0,1};\mathcal{C}_1\}$.
Straightforward evaluation of the Poisson brackets, followed by the use of Eq.~(\ref{runity}) and standard trigonometric reduction, yields
\begin{eqnarray} \nonumber
\widetilde{\mathcal{H}}'_{0,2} &=&
\mathcal{H}_{0,0}\frac{\alpha^4}{p^4}\frac{a^2}{r^2}\frac{3}{64}\frac{\eta^2}{1+\eta}  
\sum _{i=0}^2s^{2i}\sum_{j=(-1)^ii}^{i+4}\sum _{k=0}^{3-|2i-j|}B_{i,j,k}\eta^ke^{|j-2i|}\cos(jf+2ig)
\\ && \nonumber
+\mathcal{H}_{0,0}\frac{\alpha^4}{p^4}\frac{3}{8}\eta\bigg[\eta(3s^2-2)^2
+3(5s^2-4)s^2\sum_{j=1}^3\frac{2-j^\star}{j}e^{j^\star}\cos(jf+2g)\bigg]
\\ && \label{Ktidos}
+\mathcal{H}_{0,0}\frac{\alpha^4}{p^4}\frac{9}{8}\left(5 s^2-4\right)s^2
\frac{p^2}{r^2}\frac{\phi}{\eta^2}\sum _{j=1}^3(2-j^\star)e^{j^\star}\sin(jf+2g),
\end{eqnarray}
where the needed coefficients $B_{i,j,k}(s)$ are listed in Table~\ref{t:Ktidospoly}.
\par

\begin{table}[htb] 
\caption{%
Inclination polynomials $B_{i,j,k}$ in Eq.~(\protect\ref{Ktidos}). 
}
\label{t:Ktidospoly}
\small
\begin{tabular}{@{}lll@{}}
\hline\noalign{\smallskip}
$B_{0,0,0}=B_{0,0,1}$
 & $B_{1,1,0}=3 \left(37 s^2-38\right)$ & $B_{2,2,0}=B_{2,2,1}$ \\
$B_{0,0,1}=10 \left(7 s^4-16 s^2+8\right)$ & $B_{1,1,1}=12 \left(7 s^2-8\right)$ & $B_{2,2,1}=5$ \\
$B_{0,0,2}=B_{0,0,3}$ & $B_{1,1,2}=B_{1,5,0}$ & $B_{2,3,0}=B_{2,3,1}$ \\
$B_{0,0,3}=2\left(5 s^4+8 s^2-8\right)$ & $B_{1,2,0}=B_{1,2,1}$ & $B_{2,3,1}=6$ \\
$B_{0,1,0}=2 \left(57 s^4-124 s^2+60\right)$ & $B_{1,2,1}=-16 \left(4 s^2-1\right)$
  & $B_{2,3,2}=B_{2,5,2}=0$ \\
$B_{0,1,1}=4 \left(15 s^4-44 s^2+24\right)$ & $B_{1,2,2}=B_{1,2,3}$ & $B_{2,4,0}=B_{2,4,1}$ \\
$B_{0,1,2}=B_{0,3,0}$ & $B_{1,2,3}=-8 \left(s^2-2\right)$ & $B_{2,4,1}=-6$ \\
$B_{0,2,0}=2 \left(31 s^4-56 s^2+24\right)$ & $B_{1,3,0}=150-221 s^2$ & $B_{2,4,2}=B_{2,4,3}$ \\
$B_{0,2,1}=-2 \left(5 s^4+8s^2-8\right)$ & $B_{1,3,1}=-4 \left(35 s^2-24\right)$ & $B_{2,4,3}=-2$ \\
$B_{0,3,0}=2 \left(3 s^2-2\right)^2$ & $B_{1,3,2}=2-3 s^2$ & $B_{2,5,0}=B_{2,5,1}$ \\
$B_{1,-1,0}=B_{1,3,2}$ & $B_{1,4,0}=-4 \left(31 s^2-22\right)$ & $B_{2,5,1}=-10$ \\
$B_{1,0,0}=B_{1,0,1}$ & $B_{1,4,1}=-4 \left(13 s^2-10\right)$ & $B_{2,6,0}=B_{2,6,1}$ \\
$B_{1,0,1}=4 \left(15 s^2-14\right)$ & $B_{1,5,0}=-5 \left(3 s^2-2\right)$ & $B_{2,6,1}=-3$ \\
\noalign{\smallskip}\hline
\end{tabular}
\end{table}

We intentionally split $\widetilde{\mathcal{H}}'_{0,2}$ into three different blocks. Namely, all the terms on the first row of Eq.~(\ref{Ktidos}) are free from the equation of the center and factored by $a^2/r^2$, hence being of trivial integration in the true anomaly. Terms of the second row are free from both $\phi$ and $r$; the integration of terms of this type reduces to the well-known case of the integration of cosine functions in elliptic motion \cite{Brouwer1959,Kozai1962AJ,Tisserand1889}. Terms on the third row are of the form $(p/r)^2\phi\sin(mf+\alpha)$, with $m$ integer, and are readily integrated by parts. That is, on account of $\mathrm{d}(\cos{mf})/\mathrm{d}\ell=-(m/\eta^3)(p/r)^2\sin{mf}$,
\begin{equation} \label{parts}
\frac{m}{\eta^3}\int\frac{p^2}{r^2}\phi\sin{mf}\,\mathrm{d}\ell 
=-\phi\cos{mf}+\frac{\sin{mf}}{m}-\int\cos{mf}\,\mathrm{d}\ell,
\end{equation}
in this way leading to the previous case of integration of cosine functions. Particularization for definite integration follows from the fundamental theorem of calculus.
\par

It is worth noting that if, on the contrary, we replace $r$ by the conic equation in the third row of Eq.~(\ref{Ktidos}), in order to arrange this row like the Fourier series
\[
\mathcal{H}_{0,0}\frac{\alpha^4}{p^4}\frac{9}{8}(5s^2-4)s^2
\phi\sum _{j=-1}^5q_{|j-2|}e^{|j-2|}\sin(jf+2g),
\]
with $q_0=3e^2+2$, $q_1=\frac{3}{4}(e^2+4)$, $q_2=\frac{3}{2}$, and $q_3=\frac{1}{4}$, then the equation of the center shows as an isolated function of the mean anomaly when the summation index takes the value $j=0$. This arrangement brings no problem in the computation of definite integrals, which can be carried out using the general rules for computing $\langle\phi\sin{mf}\rangle_\ell$ and $\langle\phi\cos{mf}\rangle_\ell$ provided in \cite{Metris1991}. On the contrary, while indefinite integration is still possible, it requires the sophisticated use of special functions, which could make notably difficult to progress in the perturbation approach \cite{OsacarPalacian1994}.
\par

On the other hand, the evaluation of the Poisson brackets involving the integration constant $\mathcal{C}_1$ yields
\begin{eqnarray} \nonumber
\widetilde{\mathcal{H}}^*_{0,2} &=& \mathcal{H}_{0,0}\frac{R_\oplus^2}{p^2}\frac{1}{L}\frac{\partial\mathcal{C}_1}{\partial{g}}\frac{3}{2}\bigg[4-5s^2+\frac{a^2\eta}{r^2}
\sum_{i=0}^1\sum_{j=-i}^{2i+3}\sum_{k=0}^{j^\star}b_{i,j,k}\eta^{2k}e^{j'}\cos(jf+2ig)\bigg]
\\ \nonumber
&& -\mathcal{H}_{0,0}\frac{R_\oplus^2}{p^2}\frac{\partial\mathcal{C}_1}{\partial{G}}\frac{3}{2}\eta{s}^2\frac{a^2\eta}{r^2}
\sum_{j=1}^{3}[1+(j+1)^\star]e^{j^\star}\sin(jf+2g)
\\ \label{Kti2C1}
&& +\mathcal{H}_{0,0}\frac{R_\oplus^2}{p^2}\frac{\partial\mathcal{C}_1}{\partial{L}}\frac{3}{16\eta^2}\frac{a^2\eta}{r^2}
\sum_{i=0}^1b_i\sum_{j=-i}^{2i+3}q_{i,j}e^{j^\star}\sin(jf+2ig),
\end{eqnarray}
where $j'=(|j-2i|-2)j^\star$, $(j+1)^\star\equiv(j+1)\bmod2$, and the inclination polynomials $b_{i,j,k}(s)$ and the eccentricity polynomials $q_{i,j}(e)$ are provided in Table~\ref{t:C1poly}.
\par

\begin{table}[htb] 
\caption{%
Non-null inclination $b_{i,j,k}$ and eccentricity polynomials $q_{i,j}$ in Eq.~(\protect\ref{Kti2C1}). 
}
\label{t:C1poly}
\small
\begin{tabular}{llll} 
\hline\noalign{\smallskip}
$b_{0,0,0}=4-5s^2$ & $b_{1,-1,0}=\frac{1}{8}s^2$
 & \qquad & $q_{0,1}=e^2+4$ \\
$b_{0,1,0}=-\frac{1}{4}(29s^2-22)$ & $b_{1,1,0}=\frac{11}{8}s^2-1$
 & \qquad & $q_{0,2}=4 e^2$ \\
$b_{0,1,1}=\frac{1}{4}(17s^2-14)$ & $b_{1,1,1}=1-\frac{15}{8}s^2$
 & \qquad & $q_{0,3}=q_{1,-1}=e^2$ \\
$b_{0,2,0}=2-3 s^2$ & $b_{1,2,0}=5 s^2-2$
 & \qquad & $q_{1,1}=-3q_{0,1}$ \\
$b_{0,3,0}=-\frac{1}{4}(3s^2-2)$ & $b_{1,3,0}=\frac{47}{8}s^2-1$
 & \qquad & $q_{1,2}=-8(3e^2+2)$ \\
 & $b_{1,3,1}=1-\frac{19}{8}s^2$ & \qquad & $q_{1,3}=-9 q_{0,1}$ \\
 & $b_{1,4,0}=3 s^2$ & \qquad & $q_{1,4}=-24 e^2$ \\
 & $b_{1,5,0}=\frac{5}{8}s^2$ & \qquad & $q_{1,5}=-5 e^2$ \\
\noalign{\smallskip}\hline
\end{tabular}
\end{table}

In the same way as we did in the first order, we choose $\mathcal{H}_{0,2}=\langle\widetilde{\mathcal{H}}_{0,2}\rangle_\ell=\langle\widetilde{\mathcal{H}}'_{0,2}\rangle_\ell+\langle\widetilde{\mathcal{H}}^*_{0,2}\rangle_\ell$ to guarantee that it cancels all the terms of $\widetilde{\mathcal{H}}_{0,2}$ perteining to the kernel of the Lie derivative. Firstly, we compute $\langle\widetilde{\mathcal{H}}'_{0,2}\rangle_\ell$ as follows. To average the first row of Eq.~(\ref{Ktidos}) over the mean anomaly, it is first multiplied by the factor $r^2/(a^2\eta)$ to carry out the integration in the true anomaly, and then those terms that are free from $f$, which are those with $j=0$, are selected. The term free from $f$ in the second row averages to itself while the remaining terms in this row are averaged using the rule
\begin{equation} \label{Kozai}
\frac{1}{2\pi}\int_0^{2\pi}\cos(mf+\alpha)\,\mathrm{d}\ell=\left(\frac{-e}{1+\eta}\right)^m(1+m\eta)\cos\alpha,
\end{equation}
cf.~\cite{Kozai1962AJ}. Finally, the terms on the third row of Eq.~(\ref{Ktidos}) are averaged by parts with the help of Eqs.~(\ref{parts}) and (\ref{Kozai}). We finally obtain,
\begin{eqnarray} \nonumber
\langle\widetilde{\mathcal{H}}'_{0,2}\rangle_\ell &=&
\mathcal{H}_{0,0}\frac{R_\oplus^4}{p^4}\frac{3}{32}\eta\left[5(7s^4-16s^2+8)
+\eta(6s^2-4)^2+\eta^2(5 s^4+8s^2-8)\right]
\\ && \label{H02g}
+\mathcal{H}_{0,0}\frac{R_\oplus^4}{p^4}\frac{3}{16}\eta(15s^2-14)s^2e^2\cos2g,
\end{eqnarray}
which is precisely Brouwer's second-order Hamiltonian after the elimination of short-period terms. The average of Eq.~(\ref{Kti2C1}) is readily obtained with analogous procedures, to obtain
\begin{equation} \label{C1prime}
\langle\widetilde{\mathcal{H}}^*_{0,2}\rangle_\ell=
-\mathcal{H}_{0,0}\frac{R_\oplus^2}{p^2}3(5s^2-4)\frac{1}{L}\frac{\partial\mathcal{C}_1}{\partial{g}}.
\end{equation}
Visual inspection of Eqs.~(\ref{H02g}) and (\ref{C1prime}), immediately shows that if we complete the computation of the first-order term of the generating function in Eq.~(\ref{W1}) choosing
\begin{equation} \label{C1g}
\mathcal{C}_1=G\frac{R_\oplus^2}{p^2}\frac{15s^2-14}{32(5s^2-4)}s^2e^2\sin2g,
\end{equation}
then Eq.~(\ref{C1prime}) turns into the opposite of the term in the last row of Eq.~(\ref{H02g}), the only one that depends on $g$, thus mutually canceling out. Hence
\begin{equation}
\mathcal{H}_{0,2}=\mathcal{H}_{0,0}\frac{R_\oplus^4}{p^4}\frac{3}{32}\eta
\left[5(7s^4-16s^2+8)+\eta(6s^2-4)^2+\eta^2(5 s^4+8s^2-8)\right],
\end{equation}
which is the same as the second-order term of Brouwer's Hamiltonian after the elimination of long-period terms. In this way we have achieved Brouwer's total Hamiltonian reduction of the main problem at once, with a single canonical transformation.
\par

It is not a big surprise that $\mathcal{C}_1$ is the same integration constant obtained in Alfriend and Coffee's elimination of the perigee \cite{AlfriendCoffey1984,LaraSanJuanLopezOchoa2013c} or in the author's reverse normalization of the angular momentum \cite{Lara2020}, since the motion in the orbital plane is decoupled from the rotation of that plane in each case.
\par

The computation of first-order periodic corrections is now straightforward from the simple evaluation of Poisson brackets, namely $\xi-\xi'=J_2\Delta\xi$, where $\Delta\xi\equiv\{\xi,\mathcal{W}_1\}$ and $\xi$ denotes either a canonical variable or some wanted function of the canonical variables \cite{Deprit1969}. For instance, for the first-order periodic corrections to the semi-major axis we obtain
\begin{equation} \label{Da1}
\Delta{a}=a\frac{R_\oplus^2}{p^2}\frac{1}{4\eta^2}
\sum_{i=0}^1B_i(s)\sum _{j=-i}^{3+2i}A_{i,j}(\eta)e^{|j-2i|}\cos(jf+2ig),
\end{equation}
where $A_{0,0}=10-6\eta^2-4\eta^3$, $A_{0,1}=A_{1,1}=A_{1,3}=15-3\eta^2$, $A_{0,2}=A_{1,0}=A_{1,4}=6$, $A_{0,3}=A_{1,-1}=A_{1,5}=1$, and the coefficients $B_i$ are the same as those in Eq.~(\ref{Ktilde1}). Recall that Eq.~(\ref{Da1}) must be evaluated in mean (prime) variables in the direct transformation from mean elements to osculating ones, and in original (unprimed) variables in the inverse transformation from osculating to mean variables.
\par



\section{Second-order periodic corrections} \label{}

The second-order term of the generating function is now computed making $m=2$ in Eq.~(\ref{homological}). Namely
\[
\mathcal{W}_2=\mathcal{V}_2+\mathcal{C}_2, \quad \mathrm{with} \quad
\mathcal{V}_2=\frac{1}{n}\int(\widetilde{\mathcal{H}}_{0,2}-\mathcal{H}_{0,2})\mathrm{d}\ell.
\]
The needed integrals in the computation of $\mathcal{V}_2$ are either trivial, solved with the help of Eq.~(\ref{parts}) for those terms involving the equation of the center, or using the differential relation between the mean and true anomalies for those other that are free from $\phi$ but depend on trigonometric functions of $f$. Straightforward computations yield
\begin{eqnarray} \nonumber
\mathcal{V}_2 &=& G\frac{R_\oplus^4}{p^4}\frac{3\phi}{64}\bigg[
-\eta^2(5s^4+8s^2-8)-5(7s^4-16s^2+8) \\ \nonumber
&& -(15s^2-14)e^2s^2\cos2g
+12s^2(5s^2-4)\sum_{j=1}^3\frac{2-j^\star}{j}e^{j^\star}\cos(jg+2g)\bigg] \\ \label{V2poly}
&& +G\frac{R_\oplus^4}{p^4}\frac{1}{512}\sum_{i=0}^2\sum_{j=j_\mathrm{min}}^{j_\mathrm{max}}\sum_{k=0}^3
\frac{\beta_{i,j,k}(s)\eta^ks^{2i}e^{j^\star}\sin(jf+2ig)}
{(5s^2-4)^{2-i^\star}(1+\eta)^{\lfloor\frac{1}{2}(3-i)\rfloor}},
\end{eqnarray}
where $j_\mathrm{min}=2(i+1)^\star-1$, $j_\mathrm{max}=4+i+\lfloor\frac{1}{2}(i-1)\rfloor$, and the inclination polynomials $\beta_{i,j,k}$ are listed in Table~\ref{t:V2poly}.
\par

\begin{table}[htb] 
\caption{%
Non-null inclination polynomials $\beta_{i,j,k}$ in Eq.~(\protect\ref{V2poly}). 
}
\label{t:V2poly}
\small
\begin{tabular}{@{}ll@{}} 
\hline\noalign{\smallskip}
${}_{0,1,0}:-15 (3 s^2-2) (805 s^6-2448 s^4+2400
   s^2-768)$ & ${}_{1,2,3}:12 (-25 s^4+16 s^2+4)$ \\
 ${}_{0,1,1}:-3 (3 s^2-2) (2225 s^6-8160 s^4+8928
   s^2-3072)$ & ${}_{1,3,0}:2 (1855 s^4-2700 s^2+972)$ \\
 ${}_{0,1,2}:3 (-825 s^8+3030 s^6-4064 s^4+2368 s^2-512)$ & ${}_{1,3,1}:2
   (1045 s^4-1512 s^2+540)$ \\
 ${}_{0,1,3}:3 s^2 (975 s^6-2250 s^4+1728 s^2-448)$ & ${}_{1,3,2}:-2 (3
   s^2-2) (5 s^2-6)$ \\
 ${}_{0,2,0}:-\beta_{0,2,2}$ & ${}_{1,3,3}:-2 (3 s^2-2) (15
   s^2-14)$ \\
 ${}_{0,2,1}:-\beta_{0,2,3}$ & ${}_{1,4,0}:-\beta_{1,4,2}$ \\
 ${}_{0,2,2}:6 (1925 s^8-6210 s^6+7452 s^4-3936
   s^2+768)$ & ${}_{1,4,1}:-\beta_{1,4,3}$ \\
 ${}_{0,2,3}:6 (125 s^8-930 s^6+1660 s^4-1120 s^2+256)$ & ${}_{1,4,2}:-12
   (5 s^2-4) (31 s^2-22)$ \\
 ${}_{0,3,0}:-\beta_{0,3,2}$ & ${}_{1,4,3}:-12 (5 s^2-4) (13
   s^2-10)$ \\
 ${}_{0,3,1}:-\beta_{0,3,3}$ & ${}_{1,5,0}:-\beta_{1,5,2}$ \\
 ${}_{0,3,2}:2625 s^8-7270 s^6+7408 s^4-3264 s^2+512$ & ${}_{1,5,2}:-12 (3
   s^2-2) (5 s^2-4)$ \\
 ${}_{0,3,3}:s^2 (825 s^6-1990 s^4+1616
   s^2-448)$ & ${}_{2,1,0}:-\beta_{2,1,2}$ \\
 ${}_{1,-1,0}:-\beta_{1,-1,2}$ & ${}_{2,1,2}:3(225s^4-430s^2+208)$ \\
 ${}_{1,-1,1}:-\beta_{1,-1,3}$ & ${}_{2,2,0}:-\beta_{2,2,2}$ \\
 ${}_{1,-1,2}:6 (135 s^4-232 s^2+100)$ & ${}_{2,2,2}:60 (50 s^4-87
   s^2+38)$ \\
 ${}_{1,-1,3}:6 (7 s^2-6) (15 s^2-14)$ & ${}_{2,3,0}:-20 (165
   s^4-284 s^2+122)$ \\
 ${}_{1,1,0}:-24 (495 s^4-850 s^2+364)$ & ${}_{2,3,2}:8 (75 s^4-135
   s^2+61)$ \\
 ${}_{1,1,1}:-12 (855 s^4-1502 s^2+656)$ & ${}_{2,4,0}:-180
   (s^2-1) (5 s^2-4)$ \\
 ${}_{1,1,2}:48 (5 s^2-4)$ & ${}_{2,4,2}:12 (5 s^2-4) (25
   s^2-23)$ \\
 ${}_{1,1,3}:-12 (5 s^2-4) (15 s^2-14)$ & ${}_{2,5,0}:3 (5
   s^2-4) (25 s^2-18)$ \\
 ${}_{1,2,0}:12 (-95 s^4+240 s^2-132)$ & ${}_{2,5,2}:3 (5 s^2-4)
   (15 s^2-14)$ \\
 ${}_{1,2,1}:12 (-95 s^4+240 s^2-132)$ & ${}_{2,6,0}:-\beta_{2,6,2}$ \\
 ${}_{1,2,2}:\beta_{1,2,3}$ & ${}_{2,6,2}:-6 (5 s^2-4)^2$ \\
\noalign{\smallskip}\hline
\end{tabular}
\end{table}

As before, the integration constant $\mathcal{C}_2$ will be determined by imposing to the known terms of the next order
\[
\widetilde{\mathcal{H}}_{0,3}=\{\mathcal{H}_{0,2}+\mathcal{H}_{1,1},\mathcal{W}_1\}+\{\mathcal{H}_{0,1}+2\mathcal{H}_{1,0},\mathcal{W}_2\},
\]
where $\mathcal{H}_{1,1}=\mathcal{H}_{0,2}+\{\mathcal{H}_{0,1},\mathcal{W}_1\}$, the condition of being free from pure long-periodic terms. Again, the known terms are split into terms directly computable and those depending on the arbitrary function $\mathcal{C}_2$. That is, $\widetilde{\mathcal{H}}_{0,3}=\widetilde{\mathcal{H}}'_{0,3}+\widetilde{\mathcal{H}}^*_{0,3}$, where
\begin{eqnarray*}
\widetilde{\mathcal{H}}'_{0,3} &=& \{\mathcal{H}_{0,2}+\mathcal{H}_{1,1},\mathcal{W}_1\}
+\{\mathcal{H}_{0,1}+2\mathcal{H}_{1,0},\mathcal{V}_2\}, \\
\widetilde{\mathcal{H}}^*_{0,3} &=& \{\mathcal{H}_{0,1}+2\mathcal{H}_{1,0},\mathcal{C}_2\}.
\end{eqnarray*}
It follows the customary computation of $\mathcal{H}_{0,3}$ so that it cancels the terms of $\widetilde{\mathcal{H}}_{0,3}$ pertaining to the kernel of the Lie operator; namely
\[
\mathcal{H}_{0,3}=\langle\widetilde{\mathcal{H}}_{0,3}\rangle_\ell=\langle\widetilde{\mathcal{H}}'_{0,3}\rangle_\ell+\langle\widetilde{\mathcal{H}}^*_{0,3}\rangle_\ell.
\]
After straightforward evaluation of the Poisson brackets, we obtain 
\begin{equation} \label{H03l}
\langle\widetilde{\mathcal{H}}'_{0,3}\rangle_\ell=\mathcal{H}_{0,0}\frac{R_\oplus^6}{p^6}\frac{9}{512}\eta\sum_{i=0}^{2}\sum_{k=0}^{4-2i+i^\star}\frac{\beta_{i,k}(s)\eta^ks^{2i}e^{2i}}{(5s^2-4)^{2-i^\star}(1+\eta)^{i^\star}}\cos2ig,
\end{equation}
where $i^\star=i\bmod2$ and the inclination polynomials $\beta_{i,k}$ are given in Table~\ref{t:H03poly}. Analogously,
\begin{equation} \label{C2p}
\langle\widetilde{\mathcal{H}}^*_{0,3}\rangle_\ell=-\mathcal{H}_{0,0}\frac{R_\oplus^2}{p^2}\frac{9}{2}(5 s^2-4)\frac{1}{L}\frac{\partial\mathcal{C}_2}{\partial{g}}.
\end{equation}
In this process, we only found integrals of the same type as we did at the second order, and hence there were no special difficulties in solving them, yet in this case we needed to deal with notably longer series than in previous orders.
\par

\begin{table}[htb] 
\caption{%
Inclination polynomials $\beta_{i,k}$ in Eq.~(\protect\ref{H03l}). 
}
\label{t:H03poly}
\small
\begin{tabular}{@{}ll@{}} 
\hline\noalign{\smallskip}
$\beta_{0,0}=-5 \left(28700 s^{10}-107205 s^8+158960 s^6-118492 s^4+45152
   s^2-7168\right)$ \\
$\beta_{0,1}=-60 \left(3 s^2-2\right) \left(5 s^2-4\right)^2 \left(7 s^4-16 s^2+8\right)
  $ \\
$\beta_{0,2}=2 \left(28675 s^{10}-98005 s^8+130852 s^6-87164 s^4+30176 s^2-4608\right)
  $ \\
$\beta_{0,3}=-20 \left(3 s^2-2\right) \left(5 s^2-4\right)^2 \left(5 s^4+8 s^2-8\right)
  $ \\
$\beta_{0,4}=s^2 \left(15 s^2-14\right) \left(450 s^6-925 s^4+590 s^2-112\right)$ \\
$\beta_{1,0}=525 s^6-3930 s^4+5632 s^2-2256$ \\
$\beta_{1,1}=5925 s^6-16170 s^4+14848 s^2-4560$ \\
$\beta_{1,2}=\left(14-15 s^2\right) \left(75 s^4-212 s^2+120\right)$ \\
$\beta_{1,3}=\left(15 s^2-14\right) \left(45 s^4+36 s^2-56\right)$ \\
$\beta_{2,0}=\left(15 s^2-14\right)^2 \left(15 s^2-13\right)$ \\
\noalign{\smallskip}\hline
\end{tabular}
\end{table}

Once more, the simple inspection of Eqs.~(\ref{H03l}) and (\ref{C2p}) shows that if we now choose
\begin{equation} \label{C2g}
\mathcal{C}_2=G\frac{R_\oplus^4}{p^4}\frac{1}{256}\sum_{i=1}^2\sum_{k=0}^{4-2i+i^\star}\frac{\beta_{i,k}(s)\eta^ks^{2i}e^{2i}}{(5s^2-4)^{i+1}(1+\eta)^{i^\star}}\frac{\sin2ig}{2i},
\end{equation}
then Eq.~(\ref{C2p}) cancels the terms of Eq.~(\ref{H03l}) depending on the argument of the perigee out, to yield
\begin{equation}
\mathcal{H}_{0,3}=\mathcal{H}_{0,0}\frac{R_\oplus^6}{p^6}\frac{9}{512}\frac{\eta}{(5s^2-4)^{2}}\sum_{k=0}^{4}\beta_{0,k}(s)\eta^k,
\end{equation}
which is completely reduced as desired.
\par


Beyond the first order, direct and inverse transformations are no longer opposite. At the second order the direct transformation is given by $\xi=\xi'+J_2\Delta\xi+\frac{1}{2}J_2^2\delta'\xi$, where $\delta'\xi=\{\Delta\xi,\mathcal{W}_1\}+\{\xi,\mathcal{W}_2\}$ is evaluated in prime variables. The inverse transformation is $\xi'=\xi-J_2\Delta\xi+\frac{1}{2}J_2^2\delta\xi$, where $\delta\xi=\{\Delta\xi,\mathcal{W}_1\}+\{\xi,-\mathcal{W}_2\}$, is evaluated in the original variables. For instance, replacing $\xi$ by $a$ we obtain the \emph{inverse} second-order periodic correction to the semi-major axis
\begin{eqnarray} \nonumber
\delta{a} &=& a\frac{R_\oplus^4}{p^4}\frac{1}{4^4\eta^4}\bigg[
24\eta^7(5s^4+8s^2-8) +48\eta^5(15 ^2-14)s^2e^2\cos2g
\\ && \label{da2inv}
+\sum_{i=0}^2\sum_{j=-i-3i^\star}^{6+2i}\sum_{k=0}^{6-|j-2i|}
(3s^2-2)^{i^\star}s^{2i}A_{i,j,k}(s)\eta^ke^{|j-2i|}\cos(jf+2ig)
\bigg],\qquad
\end{eqnarray}
where the inclination coefficients $A_{i,j,k}$ are provided in Table~\ref{da2inv}. 
\par

\begin{table}[htb] 
\caption{%
Non-null inclination coefficients $A_{i,j,k}$ in Eq.~(\protect\ref{da2inv}). 
}
\label{t:a2poly}
\small
\begin{tabular}{@{}llll@{}} 
\hline\noalign{\smallskip}
${}_{2,-2,0}:9$ & ${}_{2,7,0}:1980$ & ${}_{1,2,0}:-11088$ & ${}_{0,0,0}:462A_{0,6,0}$ \\
${}_{2,-1,0}:108$ & ${}_{2,7,2}:-540$ & ${}_{1,2,2}:15120$ & ${}_{0,0,2}:-630A_{0,6,0}$ \\
${}_{2,0,0}:594$ & ${}_{2,8,0}:594$ & ${}_{1,2,3}:480$ & ${}_{0,0,3}:10A_{0,3,3}$ \\
${}_{2,0,2}:-54$ & ${}_{2,8,2}:-54$ & ${}_{1,2,4}:-5040$ & ${}_{0,0,4}:210A_{0,6,0}$ \\
${}_{2,1,0}:1980$ & ${}_{2,9,0}:108$ & ${}_{1,2,5}:-\frac{288(13s^2-10)}{3s^2-2}$ & ${}_{0,0,5}:24(71s^4-128s^2+56)$ \\
${}_{2,1,2}:-540$ & ${}_{2,10,0}:9$ & ${}_{1,2,6}:240$ & ${}_{0,0,6}:-4 (63 s^4-24s^2+8)$ \\
${}_{2,2,0}:4455$ & ${}_{1,-4,0}:-12$ & ${}_{1,3,0}:-9504$ & ${}_{0,1,0}:792A_{0,6,0}$ \\
${}_{2,2,2}:-2430$ & ${}_{1,-3,0}:-144$ & ${}_{1,3,2}:8640$ & ${}_{0,1,2}:-720A_{0,6,0}$ \\
${}_{2,2,4}:135$ & ${}_{1,-2,0}:-792$ & ${}_{1,3,3}:360$ & ${}_{0,1,3}:15A_{0,3,3}$ \\
${}_{2,3,0}:7128$ & ${}_{1,-2,2}:72$ & ${}_{1,3,4}:-1440$ & ${}_{0,1,4}:120A_{0,6,0}$ \\
${}_{2,3,2}:-6480$ & ${}_{1,-1,0}:-2640$ & ${}_{1,3,5}:-\frac{24(49s^2-38)}{3s^2-2}$ & ${}_{0,1,5}:-3 A_{0,3,3}$ \\
${}_{2,3,4}:1080$ & ${}_{1,-1,2}:720$ & ${}_{1,4,0}:-5940$ & ${}_{0,2,0}:495A_{0,6,0}$ \\
${}_{2,4,0}:8316$ & ${}_{1,-1,3}:24$ & ${}_{1,4,2}:3240$ & ${}_{0,2,2}:-270A_{0,6,0}$ \\
${}_{2,4,2}:-11340$ & ${}_{1,0,0}:-5940$ & ${}_{1,4,3}:144$ & ${}_{0,2,3}:6A_{0,3,3}$ \\
${}_{2,4,4}:3780$ & ${}_{1,0,2}:3240$ & ${}_{1,4,4}:-180$ & ${}_{0,2,4}:15A_{0,6,0}$ \\
${}_{2,4,6}:-180$ & ${}_{1,0,3}:144$ & ${}_{1,5,0}:-2640$ & ${}_{0,3,0}:220A_{0,6,0}$ \\
${}_{2,5,0}:7128$ & ${}_{1,0,4}:-180$ & ${}_{1,5,2}:720$ & ${}_{0,3,2}:-60A_{0,6,0}$ \\
${}_{2,5,2}:-6480$ & ${}_{1,1,0}:-9504$ & ${}_{1,5,3}:24$ & ${}_{0,3,3}:-16(3s^2-2)^2$ \\
${}_{2,5,4}:1080$ & ${}_{1,1,2}:8640$ & ${}_{1,6,0}:-792$ & ${}_{0,4,0}:66A_{0,6,0}$ \\
${}_{2,6,0}:4455$ & ${}_{1,1,3}:360$ & ${}_{1,6,2}:72$ & ${}_{0,4,2}:-6A_{0,6,0}$ \\
${}_{2,6,2}:-2430$ & ${}_{1,1,4}:-1440$ & ${}_{1,7,0}:-144$ & ${}_{0,5,0}:12A_{0,6,0}$ \\
${}_{2,6,4}:135$ & ${}_{1,1,5}:-\frac{72(43s^2-34)}{3s^2-2}$ & ${}_{1,8,0}:-12$ & ${}_{0,6,0}:2(27s^4-24s^2+8)$ \\
\noalign{\smallskip}\hline
\end{tabular}
\end{table}

\section{Initialization of the secular constants and performance tests} \label{}

Soon after Brouwer's solution appeared in print, different reports pointed out an apparent contradiction between the accuracy expected from the series truncation order and the comparatively large in-track errors obtained in a variety of tests against numerical integrations \cite{BonavitoWatsonWalden1969}. The issue, however, did not happen when fitting Brouwer's solution to observational data. Hence, the apparent discrepancy was easily identified with an inconsistency in the use of Brouwer's theory. Indeed, to get the expected accuracy provided by the secular terms, the initialization of the constants of Brouwer's solution should be done with analogous accuracy. However, Brouwer only provided the periodic corrections up to the first order of $J_2$, and hence the direct initialization of the secular mean motion for given initial conditions yields analogous accuracy. The trouble is, of course, solved if the inverse periodic corrections are computed up to the same order as the secular terms. 
\par

On the other hand, since the trouble arises from an inaccurate computation of the secular mean motion, the theory can be patched by supplementing Brouwer's first order corrections only with the inverse second-order correction to the semi-major axis, either using Eq.~(\ref{da2inv}) or in the different much shorter formulation given by \cite{LyddaneCohen1962}. Alternatively, the errors in the initialization procedure can be palliated by fitting the secular frequencies to data obtained from a preliminary numerical integration over several revolutions, or by a calibration of the secular mean motion $n'=\mu^2/L'^3$ from the energy equation \cite{BreakwellVagners1970}.
\par

The latter approach is particularly appealing because it totally avoids the need of carrying out additional computations to those already carried out by Brouwer. Thus, for given initial conditions $(\ell_0,g_0,h_0,L_0,G_0,H_0)$, the initial Hamiltonian in osculating elements evaluates to $\mathcal{H}(\ell_0,g_0,-,L_0,G_0,H_0)=E_0$. On the other hand, after the complete Hamiltonian reduction
\begin{equation} \label{energy}
E_0=-\frac{\mu^2}{2L'^2}+\sum_{m=1}^k\frac{J_2^m}{m!}\mathcal{H}_{0,m}(L',G',H)+\mathcal{O}(J_2^{k+1}).
\end{equation}
However, the constants $L'$ and $G'$ are computed from the osculating initial conditions through the inverse periodic corrections only up to $\mathcal{O}(J_2^{k-1})$. While this fact does not compromise the accuracy of Eq.~(\ref{energy}) in what respects to the terms $\mathcal{H}_{0,m}$ ($m\ge1$) because they are multiplied by corresponding factors $J_2^m$, it certainly does in the case of the Keplerian term. What Breakwell and Vagners \cite{BreakwellVagners1970} propose is then to replace $L'$ by the calibrated value
\begin{equation} \label{calibrated}
\widehat{L}=\frac{\mu}{\sqrt{2}\left[-E_0+\sum_{m=1}^k(J_2^m/m!)\mathcal{H}_{0,m}(L',G',H)\right]^{1/2}},
\end{equation}
obtained by solving the Keplerian term from Eq.~(\ref{energy}). If now $L'$ is replaced in Eq.~(\ref{energy}) by the calibrated value $\widehat{L}$ then the energy equation will remain certainly accurate to $\mathcal{O}(J_2^{k+1})$. Therefore, the initialization of the secular frequencies is notably improved using the values
\[
n_\ell=\frac{\mu^2}{\widehat{L}^3}+\sum_{m=1}^k\frac{J_2^m}{m!}\frac{\partial\mathcal{H}_{0,m}}{\partial{L'}}, \quad
n_g=\sum_{m=1}^k\frac{J_2^m}{m!}\frac{\partial\mathcal{H}_{0,m}}{\partial{G'}}, \quad
n_h=\sum_{m=1}^k\frac{J_2^m}{m!}\frac{\partial\mathcal{H}_{0,m}}{\partial{H}}.
\]
\par

Obviously, the use of Breakwell and Vagners' calibration procedure is not constrained to the second-order of Brouwer's theory, and also applies to any truncation order. In our particular case, in which we had already computed the second-order direct and inverse periodic corrections of Brouwer's theory, the calibration of the (mean) Delaunay action allowed us to improve the accuracy of Brouwer's secular terms to the third order of $J_2$ without need of computing the long series that comprise the third-order term of the generating function.
\par

We checked that the new extended third-order Brouwer's solution, in which the second-order corrections consist of a single transformation, enjoys the same accuracy as a third-order solution computed in the traditional way of splitting the Hamiltonian reduction into the sequence provided by the elimination of the parallax, followed by the elimination of the perigee, and ending with a Delaunay normalization \cite{AlfriendCoffey1984}. At the third order, the initialization of the constants of the later was also calibrated in Breakwell and Vagners' style. In both cases, the required direct and inverse transformations were computed in polar variables to avoid singularities in the case of circular orbits, but also for efficiency reasons \cite{Izsak1963AJ,Aksnes1972,Lara2015MPE,Lara2015ASR}.
\par

\begin{figure}[htb]
\begin{center}
\includegraphics[scale=1]{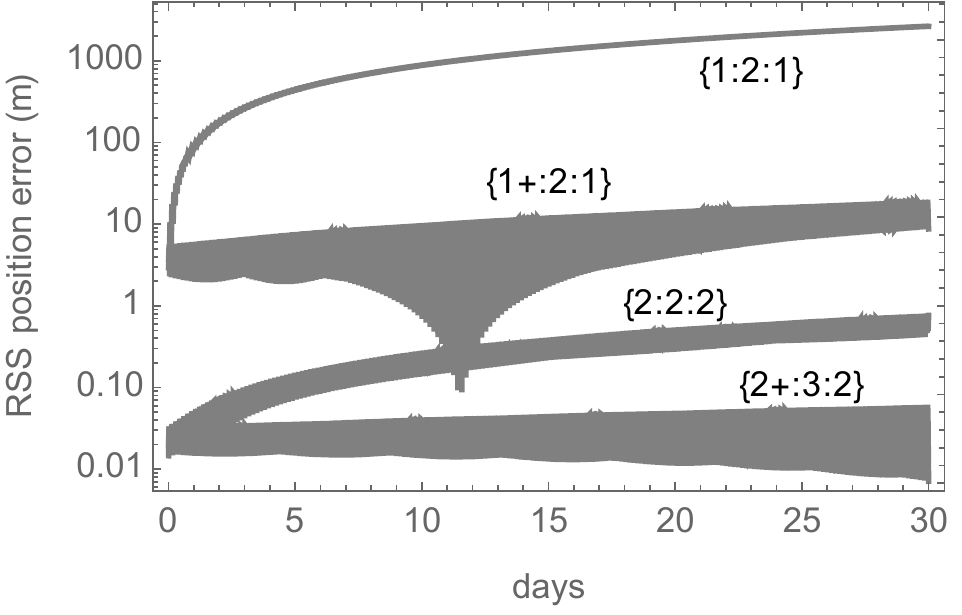}
\caption{Root-sum-square error of the Cartesian coordinates provided by different versions of Brouwer's solution. Abscissas are days.}
\label{f:errSec0xyz}
\end{center}
\end{figure}

An example of the accuracy obtained with the new single-transformation approach is presented in Fig.~\ref{f:errSec0xyz} for a Topex-type orbit, with $a=7707.270$ km, $e=0.0001$, $I=66.04$, $\Omega=180.001^\circ$, $\omega=270^\circ$, and $\ell_0=180^\circ$, yet a variety of cases have been tested for different types of orbits with analogous results. The figure depicts the root sum square (RSS) of the position errors of the analytical solution for different truncations of the extended Brouwer's solution when compared with the ``true'' solution along one month propagation. To guarantee the accuracy of the latter, the true solution was obtained from the numerical integration of the differential equations of the main problem in Cartesian coordinates using extended precision.
\par

Labels \{I:S:D\} in the plot denote the truncation order of the \textit{I}nverse corrections, \textit{S}ecular terms, and \textit{D}irect corrections of the orbital theory. The notation \{I+:S:D\} means that the inverse corrections are improved in Breakwell and Vagners' style. Thus, the label \{1:2:1\} denotes the original Brouwer's solution (with a single transformation), which at the end of one month accumulates a RSS error of about 2.5 km. The simple calibration of the secular mean motion using Eq.~(\ref{calibrated}), case \{1+:2:1\}, clearly bends the RSS errors curve towards the meter level with a negligible increase of the computational burden, reaching less than 20 m at the end of the propagation period. Figures are further improved when the orbit is propagated with the full second-order theory, in which the single second-order transformation is used both in the initialization of the constants of Brouwer's secular terms (inverse periodic corrections) and in the ephemeris computation (direct periodic corrections). Now, the RSS errors fall below the meter level at day 30, curve labeled \{2:2:2\}, yet the computational burden increases by about one third due to evaluation of second-order direct corrections, contrary to the lighter first-order corrections. Finally, supplementing Brouwer's theory with third-order secular terms and the consequent calibration of the secular mean motion, case \{2+:2:2\}, keeps the RSS errors to the level of just a few cm along the whole propagation interval with negligible increase of the computational burden with respect to the \{2:2:2\} case.
\par

For efficiency in the evaluation of perturbation solutions, arrangement of the series that comprise the periodic corrections for optimal evaluation is an important consideration \cite{CoffeyDeprit1980,HealyTravisano1998}. In this task, we limited our efforts to minor arrangements of the code, like the factorization of the inclination polynomials involved in the different summations and the following use of Horner's algorithm, and left the code optimization job to the compiler. Because we did the same for both analytical solutions (Brouwer's with single periodic corrections, and the traditional parallax-perigee-Delaunay solution), even if optimal evaluation is not achieved, the comparisons are not expected to be biased towards a particular theory.
\par

After repeated evaluation of the periodic corrections for a variety of initial conditions, we found that the evaluation of the periodic corrections of the traditional analytical solution spends roughly twice the time needed by the single-transformation approach in the evaluation of the periodic corrections. This result was a priori unexpected because the series comprising the corrections of the new approach, which only involves a single transformation, are clearly longer than the composition of those involved in the three transformations needed in the traditional approach. The improved evaluation efficiency is then attributed to the fact that the compiler is able to carry out a better optimization of the code in the case of the single-transformation approach. This fact may be understood when taking into account that, for instance, the coefficient $(5s^2-4)$ appears about 300 times in our arrangement of the periodic corrections of the single transformation, but only 73 times in the classical parallax-perigee-Delaunay transformation arranged with the same factorization criterion, where, in particular, this factor only appears in the corrections related to the elimination of the perigee. Thus, cancelling this common factor by the compiler is roughly four times more efficient in the first case. Undeniably, making a smarter organization of the code before sending it to the compiler might modify these figures. However, the balance is so radical on the side of the single transformation that these presumed improvements due to an additional preprocessing of the code are not expected to be relevant enough to revert the figures.
\par

\section{Conclusions}

Experience gained through the use of Hamiltonian simplification methods prompted us to question Brouwer's splitting normalization strategy. The convenience of dividing a normalization process into different stages has been taken for granted since the initial efforts in fully automatizing the computation of perturbation theories. Needless to say that we agree in which this way of proceeding may ease the construction of the perturbation solution. However, what is not so obvious is that the evaluation of the solution constructed this way must necessarily yield the less computational burden. On the contrary, results in this paper seem to point in the direction that the claimed benefits of partial normalization as well as Hamiltonian simplification procedures can be counterbalanced by other type of considerations, at least for the lower orders of normalization that suffice in many practical cases. Prospective application of the strategy proposed here to other instances of perturbed Keplerian motion, or to the computation of higher orders of the main problem of artificial satellite theory, should contribute to make clear the issue.
\par

Brouwer's closed-form approach and full automatization of the computation of perturbation theories seem two legitimate aims in this epoch of computational plenty. However, as demonstrated by the equation of the center controversy, rather than running perturbation algorithms in batch processes, one should not disregard the power of modern hand computations carried out with the help of existing software tools. Indeed, as far as mathematical simplification remains in the category of arts, inspection of intermediate expressions turns into a convenient practice that may eventually lead the user to straightforward simplifications that make feasible or just simpler the next step of a partially automated procedure. Like chess players, celestial mechanicians are rarely able to anticipate more than a few moves in the outcome of a perturbation approach. On the contrary, they need to wait for the opponent's reaction in order to implement a winning strategy, which, in addition, is most times settled on an empirical basis. It was, in particular, the case of the current research, in which the help provided by the computer algebra system converted into a simple task the critical inspection of the  seminal solution obtained by Brouwer.


\subsection*{Acknowledgements}

Partial support by the Spanish State Research Agency and the European Regional Development Fund under Projects ESP2016-76585-R and ESP2017-87271-P (AEI/ ERDF, EU) is recognized. The author acknowledges with pleasure the help of Sylvio Ferraz-Mello in finding particular passages of volume 2 of Poincar\'e's \textit{M\'ethodes Nouvelles}.

\end{document}